\documentclass[preprint,12pt]{elsarticle}

\usepackage{epsfig, amsmath, amssymb}

\newcounter{theoremcounter}

\newcounter{remarkcounter}
\newtheorem{theorem}[theoremcounter]{Theorem}

\newtheorem{remark}[remarkcounter]{Remark}

\newcommand{\vp}{{\mathbf p}}
\newcommand{\vz}{{\mathbf z}}
\newcommand{\vf}{{\mathbf f}}
\newcommand{\vx}{{\mathbf x}}

\def\E{\hbox{\sf E}}
\journal{Operations Research Letters}

\title{On the ergodicity bounds for  a constant retrial rate queueing   model}
\author{ A. I. Zeifman\thanks{Vologda State
University; Institute of Informatics Problems, FRC IC RAS; ISEDT RAS;
a$\_$zeifman@mail.ru}, Ya. Satin\thanks{Vologda State
University}, E.Morozov\thanks{Institute of Applied mathematical research
Karelian Research Centre RAS Center and Petrozavodsk State University}, R.Nekrasova\thanks{Institute of Applied mathematical research
Karelian Research Centre RAS Center and Petrozavodsk State University},
A.Gorshenin\thanks{Institute of Informatics Problems, FRC IC RAS}}

\date{}

\begin{document}

\maketitle

{\small

{\bf Abstract.} We consider a Markovian single-server   retrial queueing  system
with a constant retrial rate. Conditions of null ergodicity and
exponential ergodicity for the corresponding process, as well as
bounds on the rate of convergence are obtained.

\smallskip

{\bf Key words}: single-server   retrial queueing  system; constant retrial rate;
 continuous-time Markov chains;  ergodicity bounds

}

\section{Introduction}
We consider the following   Markovian single-server   retrial
queueing  system with a constant retrial rate, denoted by $\Sigma$.
The exogenous  (primary) customers follow a Poisson input with rate
$\lambda$. The customers have i.i.d. exponential service times
$\{S_i\}$, with a generic element $S$  and rate $\mu:=1/\E S $. If a
new customer finds a  server busy it joins an infinite-capacity
\textit{orbit} and called  secondary customer. We assume that the
orbit works as a single FIFS (first-in, first-served or first-come,
first-served) server. It means that, if the orbit is not empty, a
head line (the oldest) secondary customer attempts  to enter the
server after an exponentially distributed time with a rate $\mu_0$.
Thus, unlike classical retrial models, the orbit rate in $\Sigma$
does not depend on the orbit size, i.e., the number of secondary
customers. Such a model is called
 a retrial model with {\it constant retrial
rate}. Thus, orbit can be interpreted as a single-server
$\cdot/M/1$-type queue with service rate $\mu_0$, and input is
generated by the customers
 rejected in busy server. Note that the only possible source of instability of such
system is an infinite growth of the orbit size. A sufficient
stability condition of the  general single-class retrial system with
constant retrial rate described above is obtained in \cite{AvrMor}.
Stability analysis of multi-class retrial system is presented in
\cite{Seville}.

Considered system can be  successfully applied to  model the
standard multi-access protocol ALOHA, with  restrictions  for the
individual retrial rates and  has other several applications. We
list the most important papers related to motivation of presented
model. In \cite{F86} Fayolle first used  a retrial queue with
constant retrial rate to  simulate a telephone exchange system.
Then, in \cite{CSA92} the authors have modelled an unslotted Carrier
Sense Multiple Access with Collision Detection (CSMA/CD) protocol
and in \cite{CPP93} and \cite{CRP93} the authors have modelled some
particular versions of the ALOHA protocol. The model of  Fayolle
\cite{F86} was extended in \cite{A96}, \cite{AGN01} for more complex
settings, such as multiple servers and waiting places. In \cite{L96}
it has been proposed to use the retrial queue with constant retrial
rate to model a logistic system. In \cite{AY08} and \cite{AY10} the
authors have suggested to use retrial queues and retrial networks
with constant retrial rates to model TCP traffic originated from
short HTTP connections.  The present retrial model also appears to
be relevant for the optical-electrical hybrid contention resolution
scheme for Optical Packet Switching  networks
\cite{Wongetal09,Yaoetal02}.

In this note we firstly consider the description of the model
(Section 2). In Section 3, we present the auxiliary results that
were obtained for birth-death processes in \cite{gz04,z91,z06}. In
Sections 4, 5 we apply this approach for the considered model
 and obtain explicit bounds on the rate of convergence both
for null ergodic and exponentially ergodic situations.

\section{Description of the model}

Now we describe  the model in more detail. For a moment $t$, let
$\nu(t)$ be  the number of customers  in the server and  $N(t)$ be
the number of customers in the orbit. That is,   $\nu(t)=0$, if the
server is empty, and $\nu(t)=1$, otherwise, while $N(t)=0,\,
1,\,\dots$. We introduce the basic two-dimensional process $X(t)=
\{\nu(t),\,N(t),\,t\ge 0\}$  with the  state space  $\{0,\,1\}\times
\{0,\,1,\,2,\,\dots\}$. Consider in more detail the transitions
between the states of the system. To this end, we first enumerate
the states of the process as follows: each   state $\{0,\,n\}$ will
be denoted  $2n+1$, for $n\ge 0$, while  each state $\{1,\,n\}$ is
denoted $2n,\,n\ge1$. Denote by $Q=(q_{ij})$ an intensity matrix
corresponding to   the given enumeration. Then it follows that
%the $1$st line is
\begin{eqnarray*}
q_{1,1}&=&-\lambda,\\
q_{1,2}&=&\lambda,
\end{eqnarray*}
and, for $n\ge 1$,
\begin{eqnarray*}
q_{2n,2n}&=&-(\lambda+\mu),\\
q_{2n,2n+2}&=&\lambda,\\
q_{2n,2n-1}&=&\mu,
\end{eqnarray*}
\begin{eqnarray*}
q_{2n+1,2n+1}&=&-(\lambda+\mu_0),\\
q_{2n+1,2n}&=&\mu_0,\\
q_{2n+1,2n+2}&=&\lambda.
%,\,\,n\ge 1
\end{eqnarray*}
As a result, an intensity matrix $Q=(q_{ij})$ takes the form
$$
Q = \begin{pmatrix}
-\lambda & \lambda & 0 & 0 & 0 &0&0 &0&0& \cdots  \\
\mu & -(\lambda+\mu) & 0 & \lambda &0 &0& 0 &0&0& \cdots \\
0 & \mu_0 & -(\lambda+\mu_0) & \lambda &0 &0& 0&0&0& \cdots \\
0 &  0     &\mu  & -(\lambda+\mu) & 0 & \lambda & 0&0&0 & \cdots \\
0&0&0 & \mu_0 & -(\lambda+\mu_0) & \lambda &0 &0& 0& \cdots \\
0 &0&0&  0     &\mu  & -(\lambda+\mu) & 0 & \lambda & 0 & \cdots \\
\vdots &\vdots &\vdots &\vdots &\vdots &\vdots
&\vdots&\vdots&\vdots&\vdots
%\vdots & \vdots & \ddots & \vdots \\
%a_{n1} & a_{n2} & \cdots & a_{nn}
\end{pmatrix}.
$$
Figures 1,2 illustrate these transitions by different but equivalent
ways:
 \begin{figure}[!h]
   \centering
  \includegraphics[keepaspectratio=true,width=0.3\columnwidth]{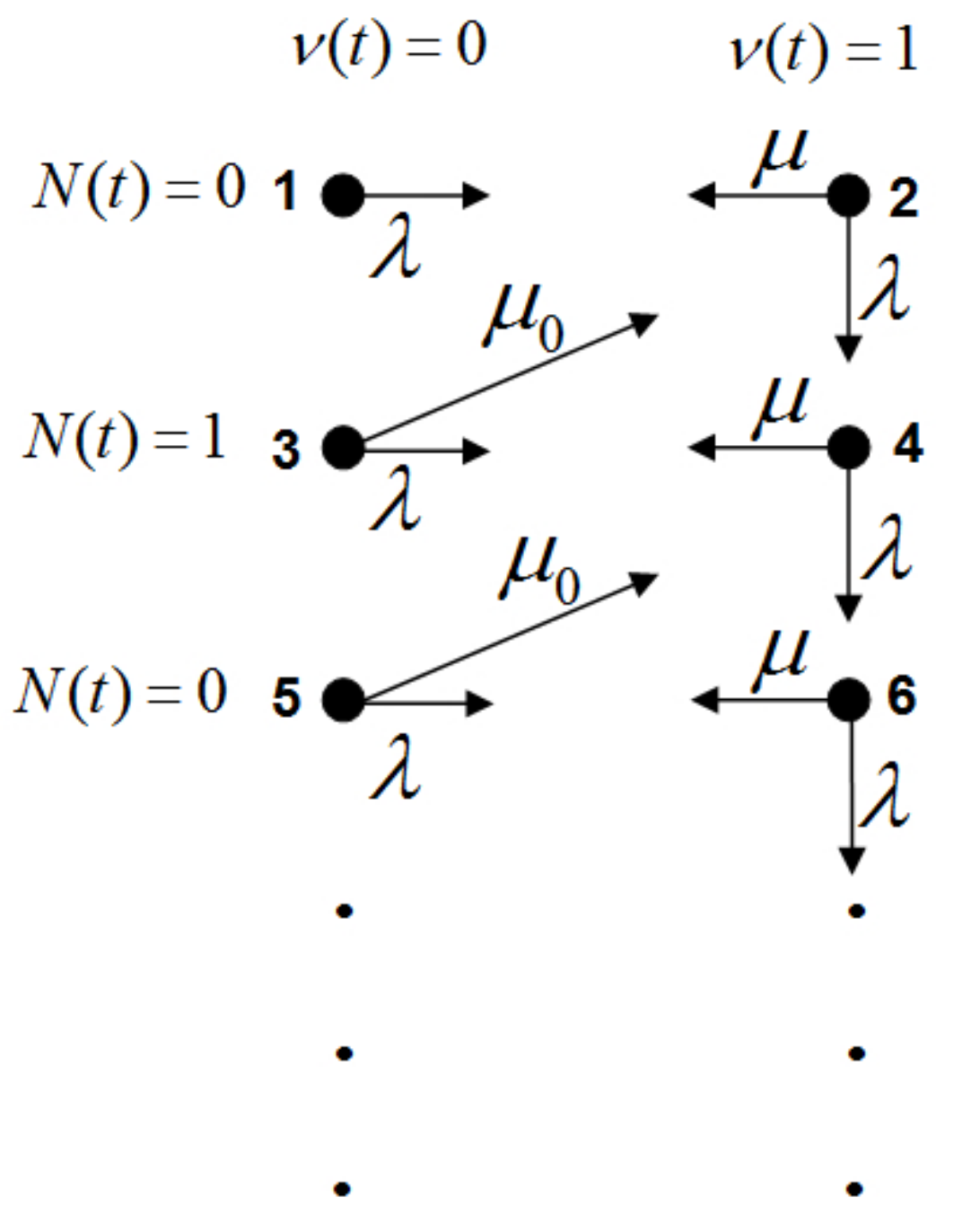}
   \caption{ }
    \label{pic1}
\end{figure}

 \begin{figure}[!h]
   \centering
  \includegraphics[keepaspectratio=true,width=0.8\columnwidth]{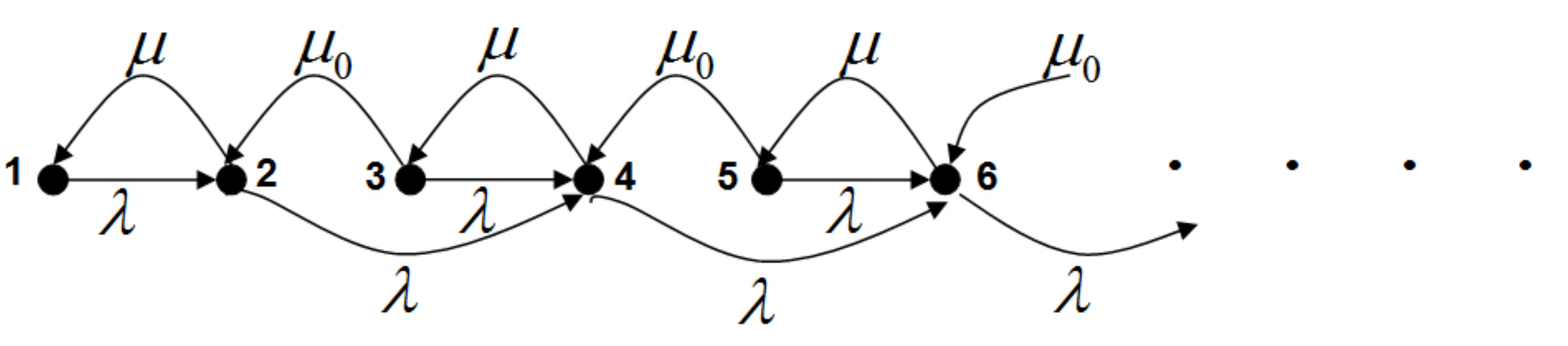}
   \caption{ }
    \label{pic1}
\end{figure}

\section{Auxiliary notions}

Then the probabilistic dynamics of the process is represented
by the forward Kolmogorov system:
\begin{equation} \label{ur01}
\frac{d\vp}{dt}=A\vp,
\end{equation}
\noindent where $A=
\left(a_{ij}\right)_{i,j=1}^{\infty}=\left(q_{ji}\right)_{i,j=1}^{\infty}
= Q^T$ is the corresponding transposed intensity matrix, and $\vp =
\vp(t) = \left(p_1(t),p_2(t),\dots\right)^T$ is the column vector of
state probabilities of the process $X(t)$.

Throughout the paper by $\|\cdot\|$  we denote  the $l_1$-norm, i.
e.  $\|{\vx}\|=\sum|x_i|$, and $\|B\| = \sup_j \sum_i |b_{ij}|$ for
matrix  $B = (b_{ij})_{i,j=1}^{\infty}$.

We recall here the approach for the study on the rate of convergence
of birth-death process, see all details and discussion in
\cite{gz04,z91,z06,z15a}.

Let $B$ be a  bounded linear operator on a Banach space ${\cal B}$
and let $I$ denote the identity operator. The number
$$
\gamma \left( B\right)_{\cal B} =\lim\limits_{h\rightarrow
+0}\frac{%
\left\| I+hB\right\| -1}h
$$
is called the logarithmic norm of $B.$

If ${\cal B} =l_1$, so that the operator $B$ is given by the matrix
$B=\left( b_{ij}\right)_{i,j=1}^\infty$, then the logarithmic norm
of $B$ can be found explicitly:
$$
\gamma \left( B\right)_1 =\sup\limits_j\left( b_{jj}
+\sum\limits_{i\neq j}\left| b_{ij} \right| \right).
$$

On the other hand, the  logarithmic norm of the operator $B$ is
related to  the Cauchy operator $V(t,s)$ of differential equation
$$\frac{d{\bf x}}{dt}=B {\bf x}$$
\noindent in the following way
$$
\gamma \left(B \right)_{\cal B} =\lim\limits_{h\rightarrow
+0}\frac{%
\left\| V\left( t+h,t\right) \right\| -1}h, \quad t\ge 0.
$$

From the latter  one can  deduce the  following bounds of the Cauchy
operator $V(t,s)$, see \cite{gz04,z91,z06}:

$$\left\| V\left( t,s\right) \right\|_{\cal B} \leq e^{(t-s)\gamma \left(
B \right)  },\quad 0\le s\le t.  $$

\medskip

Let $\Omega$ be a set of all stochastic vectors, i. e. $l_1$ vectors
with nonnegative coordinates and unit norm. Hence $\|A\| = 2\lambda
+ 2\max(\mu,\mu_0) < \infty$. Thus, the Cauchy problem for
differential equation (\ref{ur01}) has a unique solutions for an
arbitrary initial condition, and  $\vp(s) \in \Omega$ implies
$\vp(t) \in \Omega$ for $t \ge s \ge 0$.

\section{Null ergodicity}

We recall that a Markov chain $X(t)$ is called null ergodic (in
other terms, it means null recurrence or  transience), if $p_k(t)
\to 0$ as $t \to \infty$ for any initial condition ${\bf p}(0)$ and
any $k$.

\smallskip

Consider a decreasing  sequence of positive numbers $\{\delta_i\}$,
$i=1,2, \dots$,  and the corresponding diagonal matrix $\Delta$ with
diagonal entries $\{\delta_k\}$.

Let $l_{1\Delta}$ be the space of sequences:
$$
l_{1\Delta}= \left\{{\bf p}= (p_1,p_2,\cdots)^T :\, \|{\bf
p}\|_{1\Delta} \equiv \|\Delta {\bf p}\| <\infty \right\}. $$

Then we have
\begin{equation}
\gamma \left( A \right)_{1\Delta} =\gamma \left(\Delta A
\Delta^{-1}\right)_1 = \sup_i\left(a_{ii}+\sum_{j\neq
i}\frac{\delta_{j}}{\delta_i}a_{ji}\right). \label{nul01}
\end{equation}

Put $\delta_{1} = 1$ and $\delta_{2k+1}=a\delta_{2k}$,
$\delta_{2k}=b\delta_{2k-1}$, if $k \ge 1$ for some $a,b$. Then
$\delta_{2k+1}=a^kb^k$, $\delta_{2k}=a^{k-1}b^k$, if $k \ge 1$.

We have
$$
A=Q^T = \begin{pmatrix}
-\lambda & \mu & 0 & 0 & 0 &0&0 &0&0& \cdots  \\
\lambda & -(\lambda+\mu) &  \mu_0 &0 &0 &0& 0 &0&0& \cdots \\
0 & 0 & -(\lambda+\mu_0) & \mu &0 &0& 0&0&0& \cdots \\
0 &  \lambda     &\lambda  & -(\lambda+\mu) & \mu_0 & 0 & 0&0&0 & \cdots \\
0&0&0 & 0 & -(\lambda+\mu_0) & \mu &0 &0& 0& \cdots \\
0 &0&0&  \lambda     &\lambda  & -(\lambda+\mu) & \mu_0 & 0 & 0 & \cdots \\
\vdots &\vdots &\vdots &\vdots &\vdots &\vdots
&\vdots&\vdots&\vdots&\vdots
%\vdots & \vdots & \ddots & \vdots \\
%a_{n1} & a_{n2} & \cdots & a_{nn}
\end{pmatrix},
$$
\noindent and
$$
\Delta A \Delta^{-1} = \begin{pmatrix}
-\lambda & b^{-1}\mu & 0 & 0 & 0 &0&0 &0&0& \cdots  \\
b\lambda & -(\lambda+\mu) &  a^{-1}\mu_0 &0 &0 &0& 0 &0&0& \cdots \\
0 & 0 & -(\lambda+\mu_0) & b^{-1}\mu &0 &0& 0&0&0& \cdots \\
0 &  ab\lambda     &b\lambda  & -(\lambda+\mu) & a^{-1}\mu_0 & 0 & 0&0&0 & \cdots \\
0&0&0 & 0 & -(\lambda+\mu_0) & b^{-1}\mu &0 &0& 0& \cdots \\
0 &0&0&  ab\lambda     &b\lambda  & -(\lambda+\mu) & a^{-1}\mu_0 & 0 & 0 & \cdots \\
\vdots &\vdots &\vdots &\vdots &\vdots &\vdots
&\vdots&\vdots&\vdots&\vdots
%\vdots & \vdots & \ddots & \vdots \\
%a_{n1} & a_{n2} & \cdots & a_{nn}
\end{pmatrix}.
$$

Therefore
\begin{eqnarray}
\gamma \left( A \right)_{1\Delta} = - \min \left(\lambda(1-b),
\lambda(1-ab)-\mu(b^{-1}-1),\lambda(1-b)-\mu_0(a^{-1}-1)\right)=
\nonumber \\ - \min
\left(\lambda(1-ab)-\mu(b^{-1}-1),\lambda(1-b)-\mu_0(a^{-1}-1)\right).
\label{nul02}
\end{eqnarray}

\bigskip

Let now
\begin{equation}
\mu \mu_0 < \lambda(\lambda+\mu_0). \label{nul05}
\end{equation}

Then one has

\begin{equation} b^*: =
\frac{\mu(\lambda+\mu_0)}{\lambda(\lambda+\mu+\mu_0)} <1.
\label{nul04}
\end{equation}

Now for any $b \in (b^*,1)$ the inequality
\begin{equation}
\frac{\mu_0}{\lambda(1-b)+\mu_0} <
\frac{\lambda-\mu(b^{-1}-1)}{b\lambda} \label{nul031}
\end{equation}
\noindent holds.

Choose  now

\begin{equation}
a \in \left( \frac{\mu_0}{\lambda(1-b)+\mu_0},
\frac{\lambda-\mu(b^{-1}-1)}{b\lambda}\right). \label{nul03}
\end{equation}

Then $a<1,b<1$, $\lambda(1-ab)-\mu(b^{-1}-1) >0$,
$\lambda(1-b)-\mu_0(a^{-1}-1)>0$.

Hence inequality $\gamma \left( A \right)_{1\Delta} < 0$ implies
null ergodicity of  $X(t)$, and we obtain the following statement.

\begin{theorem}
Let assumption (\ref{nul05}) hold. Then the process $X(t)$ is null
ergodic and
\begin{equation}
\sum_{i=0}^{N} p_i(t) \le \frac{\delta_k}{\delta_N}e^{- \zeta^* t},
\label{nul06}
\end{equation}
\noindent for any $ t \ge 0$, any initial condition $X(0)= k$, and
any natural $N$,
where $$\zeta^*= \min \left(\lambda(1-ab)-\mu(b^{-1}-1),\lambda(1-b)-\mu_0(a^{-1}-1)\right) > 0.$$
\end{theorem}

\section{Exponential ergodicity}

By introducing $p_1(t) = 1 - \sum_{i \ge 2} p_i(t)$, from
(\ref{ur01}) we obtain the equation
\begin{equation}
\frac{d\vz}{dt}= B\vz+\vf, \label{erg01}
\end{equation}
\noindent where $\vf =\left( a_{21}, a_{31}, \ldots
\right)^{T}=\left( \lambda, 0, 0, \ldots \right)^{T} $, $\vz =
\vz(t)=\left(p_2(t), p_3(t),\cdots \right)^{T}$,
\begin{equation}
{%\scriptsize
B = \left(b_{ij}\right)_{i,j=1}^{\infty} = \left(
\begin{array}{ccccc}
a_{22} -a_{21} & a_{23} -a_{21} & \cdots & a_{2r}
-a_{21} & \cdots \\
a_{32} -a_{31} & a_{33} -a_{31} & \cdots & a_{3r}
-a_{31} & \cdots \\
\cdots & \cdots & \cdots & \cdots  & \cdots \\
a_{r2} -a_{r1} & a_{r3} -a_{r1} & \cdots & a_{rr}
-a_{r1} & \cdots \\
\cdots & \cdots & \cdots & \cdots  & \cdots
\end{array}
\right) ,}\label{erg02}
\end{equation}
see detailed discussion in \cite{z91,z06,z15a}. Let $\{d_i\}$,
$i=2,3, \dots$, be a  sequence of positive numbers. Put
 $g_i=\sum_{n=2}^i d_n.$

Let $D$ be the upper triangular matrix,
\begin{equation}
D=\left(
\begin{array}{ccccccc}
d_2   & d_2 & d_2 & \cdots  \\
0   & d_3  & d_3  &   \cdots  \\
0   & 0  & d_4  &   \cdots  \\
& \ddots & \ddots & \ddots \\
\end{array}
\right), \label{erg03}
\end{equation} \noindent  and $l_{1D}$ be the
corresponding space of sequences
$$l_{1D}=\left\{{\bf z} = (p_2,p_3,\cdots)^{T} |\, \|{\bf z}\|_{1D} \equiv \|D {\bf z}\|_1 <\infty \right\}.$$

Consider equation (\ref{erg01}) in the space $l_{1D}$. Then the
logarithmic norm $\gamma({B})_{1D} = \gamma(DBD^{-1})_1$, see
\cite{gz04,z06}.

We have
$$
D B D^{-1} = \begin{pmatrix}
-(\lambda+\mu) & \frac{d_2}{d_3}\mu & 0 & 0 & 0 &0&0 &0&0& \cdots  \\
\frac{d_3}{d_2}\lambda & -(\lambda+\mu_0) &  \frac{d_3}{d_4}\mu_0 &0 &0 &0& 0 &0&0& \cdots \\
\frac{d_4}{d_2}\lambda & 0 & -(\lambda+\mu) & \frac{d_4}{d_5}\mu &0 &0& 0&0&0& \cdots \\
0 &  0    &\frac{d_5}{d_4}\lambda  & -(\lambda+\mu_0) & \frac{d_5}{d_6}\mu_0 & 0 & 0&0&0 & \cdots \\
0&0& \frac{d_6}{d_4}\lambda &0 & -(\lambda+\mu) & \frac{d_6}{d_7}\mu &0 &0& 0& \cdots \\
\vdots &\vdots &\vdots &\vdots &\vdots &\vdots \vdots &\vdots &\vdots &\vdots \\
\vdots &\vdots &\vdots &\vdots &\vdots &\vdots
&\vdots&\vdots&\vdots&\vdots
%\vdots & \vdots & \ddots & \vdots \\
%a_{n1} & a_{n2} & \cdots & a_{nn}
\end{pmatrix}.
$$

Put $d_{2} = 1$ and $d_{2k+1}=bd_{2k}$, $d_{2k+2}=ad_{2k+1}$, if $k
\ge 1$.

Then
\begin{equation}
\gamma \left( B \right)_{1D} = - \inf_{i \ge 2} \alpha_i,
\label{erg04}
\end{equation}
\noindent where
\begin{equation}
\begin{array}{ccc}
\alpha_2 = \lambda + \mu - \lambda(b + ab),\\
\alpha_{2k+1} = \lambda + \mu_0 - \mu b^{-1}, \quad k\ge 1, \\
\alpha_{2k+2} = \lambda + \mu - \lambda(b + ab) - \mu_0 a^{-1},
\quad k\ge 1.
\end{array}
\label{erg05}
\end{equation}
Hence exponential  ergodicity of the process $X(t)$ follows from the
bound
\begin{equation}
\inf_{i \ge 1} \alpha_i = \min \left(\lambda + \mu_0 - \mu b^{-1},
\lambda + \mu - \lambda(b + ab) - \mu_0 a^{-1} \right) > 0,
\label{erg06}
\end{equation}
\noindent for some $a,b$ such that $ab > 1$, see
\cite{gz04,z91,z06}.

Let now that
\begin{equation}
\mu\mu_0 > \lambda(\lambda+\mu_0), \label{erg080}
\end{equation}
\noindent in contrast to (\ref{nul05})

 Put $x=ab$.

Then one has
\begin{equation}
x^*: = \frac{\mu\mu_0}{\lambda(\lambda+\mu_0)} >1. \label{erg08}
\end{equation}

Then  for any $x \in (1,x^*)$ the inequality
\begin{equation}
\frac{\mu}{\lambda + \mu_0} < \frac{\lambda+\mu - x \lambda}{\lambda
+\mu_0/x}. \label{erg071}
\end{equation}
\noindent holds.

Choose  now

\begin{equation}
b \in \left(\frac{\mu}{\lambda + \mu_0},  \frac{\lambda+\mu - x
\lambda}{\lambda +\mu_0/x}\right), \label{erg07}
\end{equation}
\noindent and $a = x/b$.

Then $\lambda + \mu_0 - \mu b^{-1} > 0,$ $ \lambda + \mu - \lambda(b
+ ab) - \mu_0 a^{-1} > 0$, and as a result,  we obtain the following
statement.

\begin{theorem}
Let assumption (\ref{erg080}) be true. Then the process $X(t)$ is
exponentially ergodic and the following bound holds:
\begin{equation}
\|{\bf p}(t) - {\bf \pi}\|_1  \le 4 e^{-\alpha^* t}\sum_{i \ge 2}g_i
|p_i(0) - \pi_i|, \label{erg09}
\end{equation}
\noindent for any $ t \ge 0$, and any initial condition ${\bf p}(0)$, where
$$\alpha^* = \min \left(\lambda + \mu_0 - \mu b^{-1}, \lambda + \mu - \lambda(b + ab) - \mu_0 a^{-1} \right) > 0,$$
and ${\bf \pi} = \left(\pi_0,\pi_1,\dots\right)^T$ is the
corresponding stationary distribution.
\end{theorem}

\medskip

\begin{remark} It is worth mentioning that condition (\ref{erg080}) is the
stability criteria  of the system studied by regenerative method in  \cite{AvrMor}.
More exactly, under this condition, the workload and queue size
processes are positive recurrent regenerative (that is, the mean
regeneration period is finite).
\end{remark}

\begin{remark} One can obtain the respective perturbation bounds
applying Theorem 2 and approach of \cite{z14c}.
On the other hand, the opposite condition (\ref{nul05}) implies
an unlimited growth of the process, and it is consistent with the null
ergodicity which established in Theorem 1.
\end{remark}

\section{Acknowledgement}
The research is supported by the Russian Foundation for Basic
Research, projects no. 15-01-01698,15-07-02341,15-37-20851; and by
Ministry of Education and Science, State Contract No. 1816.

\end{document}